\documentclass[11pt]{amsart}
\usepackage{amsfonts,amsmath,amssymb}
\usepackage[all]{xy}
\usepackage{hyperref}
\begin{document}

\newcommand{\ie}{\textit{i}.\textit{e}.\,}
\newcommand{\eg}{\textit{e}.\textit{g}.\,}
\newcommand{\cf}{{\textit{cf}.\,}}
\parindent=0pt 
\parskip=6pt

\newcommand{\la}{{\langle}}
\newcommand{\ra}{{\rangle}}
\newcommand{\A}{{\mathcal A}}
\newcommand{\B}{{\mathcal B}}
\newcommand{\BFK}{{\rm BFK}}
\newcommand{\bk}{{\boldsymbol{\kappa}}}
\newcommand{\br}{{\boldsymbol{\rho}}}
\newcommand{\bo}{{\boldsymbol{\omega}}}
\newcommand{\bbar}{{\bar{\omega}}}
\newcommand{\ba}{{\bf a}}
\newcommand{\bc}{{\bf c}}
\newcommand{\bi}{{\bf i}}
\newcommand{\bq}{{\bf k}}
\newcommand{\bm}{{\bf m}}
\newcommand{\bn}{{\bf n}}
\newcommand{\bu}{{\bf u}}
\newcommand{\bv}{{\bf v}}
\newcommand{\bY}{{\bf Y}}
\newcommand{\cZ}{{\mathcal Z}}
\newcommand{\C}{{\mathbb C}}
\newcommand{\cL}{{\mathcal L}}
\newcommand{\cHd}{{{\mathcal H}^{dif}}}
\newcommand{\cHi}{{{\mathcal H}^{inv}}}
\newcommand{\End}{{\rm End}}
\newcommand{\half}{{\textstyle{\frac{1}{2}}}}
\newcommand{\M}{{\mathbb M}}
\newcommand{\N}{{\mathbb N}}
\newcommand{\Q}{{\mathbb Q}}
\newcommand{\R}{{\mathbb R}}
\newcommand{\T}{{\mathbb T}}
\newcommand{\sT}{{\sf T}}
\newcommand{\Z}{{\mathbb Z}}
\newcommand{\Tor}{{\rm Tor}}
\newcommand{\Hom}{{\rm Hom}}
\newcommand{\Maps}{{\rm Maps}}
\newcommand{\Cay}{{\rm Cay}}
\newcommand{\Lie}{{\rm Lie \;}}
\newcommand{\Sp}{{\rm Sp}}
\newcommand{\U}{{\rm U}}
\newcommand{\op}{{\rm op}}
\newcommand{\dR}{{\rm dR}}
\newcommand{\res}{{\rm res}}
\newcommand{\QS}{{\sf Q}}
\newcommand{\NS}{{\sf N}}
\newcommand{\sS}{{\sf S}}
\newcommand{\iso}{{\rm iso}}
\newcommand{\SD}{{\rm SDiff}}
\newcommand{\FN}{{\bf{N}_*}}
\newcommand{\NF}{{\bf N}}

\title{Topological invariants of some chemical reaction networks }

\author[Jack Morava]{Jack Morava}

\address{Department of Mathematics, The Johns Hopkins University,
Baltimore, Maryland 21218}

\email{jack@math.jhu.edu}

    
                                                                
                                                                  

\date{June 2020}

\begin{abstract}{Certain toric dynamical systems \cite{15,21,22,25} studied
in physical chemistry have associated toric varieties which, when smooth,
represent elements in the homotopy groups $M\xi_*B\T$ of a symplectic 
variant of the $A_\infty$ Baker-Richter spectrum $M\xi$ \cite{9,32}. The 
noncommutative Hopf algebroid $(M\xi_*,M\xi_*M\xi)$ has interesting
connections to the group of formal diffeomorphism of the noncommutative 
line, conjecturally defining a noncommmutative analog of the  Helmholtz 
free energy for systems of complex symplectic (completely integrable) 
Hamiltonian toric manifolds.}\end{abstract}

\maketitle 

{\bf Introduction and organization} 

In 2006 Baker and Richter \cite{8,9} defined an $A_\infty$ Thom spectrum 
$M\xi$ associated to the loop space $\Omega \Sigma \C P^\infty$, with close 
relations to quasi- and noncommutative symmetric functions, and the 
remarkable property that it splits $p$-locally as a wedge of copies of 
$BP$, even though $M\xi$, being noncommutative, is not an $M\U$-algebra.  
Building on this work a few years later, N Kitchloo and the author \cite{32}
sketched a theory of $H_*M\xi$-valued characteristic numbers for the complex 
quasi-toric manifolds defined by Davis and Januszkiewicz \cite{16}. 

This paper applies these ideas to work of Guillemin, Ginzburg, and Karshon 
\cite{1} on cobordism of Hamiltonian group actions; in particular, to certain 
complex symplectic toric manifolds related to classical completely integrable
systems. Smooth Hamiltonian toric manifolds associated to chemical reaction 
networks \cite{15,21,22,25} provide a significant class of examples.

Part I of this paper begins with sketches (\S 1.1, 1.2) of background material 
from geometry and mechanics, using the language of algebraic topology; 
the recent monographs of Guillemin {\it et al}, and of Buchstaber and Panov 
\cite{2} on toric geometry, make this practical. A summary of some 
modern work on chemical reaction networks follows, and the final section 
reviews characteristic classes and numbers from non-equivariant cobordism 
theory. A brief account of the motivation behind this paper is postponed 
to an afterword. 

Part II is concerned more directly with the spectrum $M\xi$ and its associated
noncommutative Hopf algebroid $(M\xi_*,M\xi_*M\xi)\otimes \Q$ of rational 
cooperations, conjecturally closely related to a noncommutative Hopf algebra 
of formal diffeomorphisms of the noncommutative line \cite{10}. Its 
connections with statistical mechanics \cite{30,36,37} and free probability 
theory \cite{19,20,26,29,41} suggest a candidate for the Helmholtz free energy 
for systems of symplectic Hamiltonian toric manifolds. 

{\bf Apology} to the reader: This paper takes its target audience to be 
generally familiar with modern algebraic topology (such as cobordism theory 
and equivariant language), with an interest in applications. Some knowledge 
of toric geometry, classical mechanics, and combinatorics would also be useful
but \S 1.1-2 try to provide enough of an Ariadne's thread through current 
literature to make the rest of the paper accessible. The results of the 
paper involve the algebraic theory of chemical reaction networks and 
quasi-symmetric functions, neither of which is assumed to be known; 
\S 1.3 and 2.1 try to sketch enough of these topics for the reader to see 
why the proposed applications may be interesting.
   
{\bf Acknowledgements} The author has profited from conversations and
correspondence with many helpful colleagues over the years, and would
especially like to thank Andrew Baker, Birgit Richter, Nitu Kitchloo, 
Michiel Hazewinkel, John McKay, Roland Friedrich, Masoud Khalkali, and
Jeremy Gunawardena for their interest and help; but of course the
excesses and misunderstandings in this paper are his own responsibiity.
He moreover thanks the editors and  referees for their insightful comments, 
recommendations, and patience. \bigskip

{\bf \S I Background} 

{\bf 1.1 Some toric geometry} 

{\bf 1.1.1} This paper is concerned with applications of the theory of 
Hamiltonian $G$-cobordism developed in [1], specialized to the context of 
Delzant-like half-dimensional toric manifolds. The general theory of 
Guillemin, Ginzburg, and Karshon is in some sense naturally bigraded, by 
symmetry group as well as dimension: the product of an $n$-dimensional 
$G$-manifold with an $m$-dimensional $H$-manifold is most naturally an 
$n+m$-dimensional $G \times H$ manifold. Here our focus will be a neighborhood
of the diagonal of this bigraded category.

We will be interested in compact smooth `$\sT$-manifolds' $X$ endowed with 
an effective\begin{footnote}{\ie $(\forall x \in X) \; \Lie T \to \End(T_x X)$ 
is injective; we will also fix an isomorphism of Lie $T$ with $\pi_1(T) 
\otimes \R.$ We write $\T = \R/\Z$, $T$ for a generic torus, and $T^k \cong 
\T^k$ for a specific torus.}\end{footnote} smooth action of a compact torus 
$T$, such that $\dim T = \lfloor \half \dim X \rfloor$; thus the toruses 
acting on $\sT$-manifolds $(X_0,T_0), \; (X_1,T_1)$ will be distinct. 
The boundary $\partial X$ of an odd-dimensional $\sT$-manifold is again a 
$\sT$-manifold, and if one or both of a pair of $\sT$-manifolds as above is 
even-dimensional, then $M_0 \times M_1$ is a $\sT$-manifold with an action of 
$T_0 \times T_1$ and boundary $\partial M_0 \times M_1 \coprod M_0 \times 
\partial M_1$. The class of $\sT$-manifolds is thus a cobordism category with 
products, graded by dimension: an object is a smooth manifold with an action 
by a torus of half its dimension, much like the spaces classified by Delzant 
\cite{1}(\S 1.3 p 17, \S 1.5 p 29), while odd-dimensional $\sT$-manifolds
will be the morphisms of the category. This can be elaborated by specifying 
various orientations, but terminology for structures (symplectic, 
complex-oriented \dots) on manifolds with group actions can become confusing. 
Our convention is that a $\mathcal{Z}$-structure on a $\sT$-manifold will 
always mean a $T$-equivariant $\mathcal{Z}$-structure. 

To define the orientations relevant below, we work with $2n$-dimensional
$\sT$-manifolds with an underlying Davis-Januszkiewicz quasi-toric 
structure, \ie such that the quotient
\[
M \to M/T = P
\]
is homeomorphic to a simple (\ie with $n$ facets meeting at each vertex 
\cite{16}(\S 1.13),\cite{2}(1.1 p 4, 1.7 p 46) polyhedron, with boundary dual 
to a  simplicial polytope $K$ with a set $[m]$ of vertices. The map which 
assigns to a facet of $P$, its isotropy subgroup, lifts (nonuniquely) to an 
{\it omniorientation}: a surjective homomorphism 
\[
\Lambda :\Z^m \to \Z^n \cong \pi_1(T)
\]
\cite{2}(\S 7.3.9 p 247) such that the minor $\Lambda_{\{i_v\}}$ associated 
to the columns $\{i_v\}$ indexed by the facets meeting at the vertex $v$ has 
determinant $\pm 1$. Such an omniorientation defines, and is defined uniquely 
by, an equivariant complex orientation \cite{43}, \ie a stable complex 
structure \cite{2}(B.6.1 p 434] with an associated equivariant decomposition 
\[
T_M \oplus \C^{m-n} \cong \oplus_{1 \leq i \leq m} L_i
\]
of the tangent bundle of $M$ as a sum of complex line bundles; see \cite{7}, 
and \S 1.2.3 below. 

{\bf 1.1.2 Definition} An (integral) {\it Hamiltonian} (complex symplectic) 
orientation on $M$ is defined by 

1) an equivariant complex line bundle $\cL$ on $M$, with characteristic 
class $c^T_1(\cL) \in H^2_T(M;\Z)$ \cite{1}(C.47 p 214), together with a 
lift of $c^T_1(\cL)$ to

2) an equivariant symplectic form $\omega \in \Omega^2_T(M)$, \ie $d\omega =
0$ with $\omega^n$ everywhere nonzero, together with an equivariantly closed
2-form $\bbar = \omega - \Phi$ such that $(2\pi)^{-1}[\bbar]$ is the image
in $H^2_T(M;\R)$ of $c^T_1(\cL)$. Then 

3) $\Phi : M \to (\Lie T)^*$ is a moment map for the action of $T$ on $M$
\cite{1}(\S 2.15 p 25, E.37 p 266) identifying $P$ with a convex polytope
in $(\Lie T)^*$; moreover, there is a homotopically unique almost complex 
structure on $M$ compatible with $\omega$, defining a canonical complex 
omniorientation on $M$.

This is a generalization of the stable complex Hamiltonian $T$-structure 
on $M$ defined by a smooth projective toric variety with an explicit 
omniorientation\begin{footnote}{It is tempting to think of Hamiltonian
toric manifolds as something like compactifications of classical harmonic
oscillators, much as toric varieties are compactifications of algebraic
torus actions}\end{footnote}. A complex-oriented $2n+1$-manifold with $T^n$ 
action, quasi-toric on its boundary, together with an equivariant complex 
line bundle and compatible two-forms $(\bbar,\omega)$ as above, nondegenerate 
\cite{1}(Thm H.10 p 309) in the sense that $\omega^n$ has one-dimensional 
kernel (trivial on the boundary), defines a symplectic cobordism with a
complex-oriented version of a formal Hamiltonian structure as in 
\cite{1}(E.47 p 273, H1.1 p 301). Relaxing the integrality condition 
on $[\bbar]$ in effect tensors the underlying homotopy-theoretic 
constructions with $\R$. 

There is however an interesting alternative notion of a complex-oriented
cobordism between such manifolds, defined not by a generalized moment map
\cite{1}(3.21 p 35) but by a stable decomposition of its tangent bundle into
complex line bundles, as above. This leads to the quite different theory of
characteristic numbers associated to the ring-spectrum $M\xi$ pursued in
this paper.

{\bf 1.1.3 Examples} 

{\it i)} A projective toric manifold $M$ (defined by a smooth complex toric 
variety) inherits a K\"ahler class $\omega$ from its projective embedding 
\cite{2}(\S 5.5.4 p 197, \S 7.3 p 243), as well as a moment map $\Phi$ defining
a polyhedron in $(\Lie T)^*$ bounded by $m$ half-planes $\{\bv \:|\: \bv \cdot 
e_i \geq \lambda_i\}$, such that \cite{24}
\[
(2\pi)^{-1}[\omega] = - \sum \lambda_i [v_i] \in H^2(M,\R)
\]
for certain classes $[v_i] \in H^2(M,\Z)$ defined in \S 1.2.2 below; this is
generalized in \cite{11}. Similarly, toric orbifolds have underlying toric
varieties. Note that $(M,-\omega)$ is not projective. 

{\it ii)} The orbits of the coadjoint action of a compact Lie group $G$ on 
the dual $(\Lie \: G)^*$ of its Lie algebra are canonically symplectic 
manifolds, with the inclusion
\[
G/\iso(\lambda) \cong G \cdot \lambda \subset (\Lie G)^*
\]
defining a moment map\begin{footnote}{See \cite{49} for algebraic groups}\end
{footnote}. They play a role in the representation theory of general Lie
groups, in Kirillov's program, analogous to that of conjugacy classes in
the theory of finite groups. Such orbits are foliated by tori, and are a 
rich source of classical (Liouville - Arnol'd) completely integrable 
Hamiltonian dynamical systems, \eg the Eulerian mechanics of rigid bodies. 
This is pursued further in \S III. \bigskip

{\bf 1.2 Some toric topology} 

Relations between the toric manifolds of Davis and Januszkiewicz [16] and
the toric varieties of algebraic geometry \cite{14,48} are a source of possible
ambiguity in the discussion below; \cite{2} is helpful, especially for
connections with combinatorics. From the point of view of homotopy theory,
the relevant fact is that toric objects are abundantly supplied with 
interesting line bundles, conveniently indexed in various ways by polyopes.
The characteristic classes of these line bundles provide algebraic 
invariants defined over the integers, and thus over $\Q$; but in the latter
context distinctions between manifolds, orbifolds, and varieties \cite{23,34}
can sometimes be usefully backgrounded. \bigskip 

{\bf 1.2.1} Let $K$ be a finite simplicial complex with vertices $[m]$, 
and let $(X_i,A_i)$ be a family of reasonable pairs of spaces, also indexed 
by $[m]$; then Browder's polyhedral product 
\[
K \mapsto (X_*,A_*)^K
\]
(\cite{2}(\S 4.2.3 p 137), functorial under maps of families $(X_*,A_*)$ and 
inclusions of complexes) associates to $K$, the moment-angle complex 
\cite{2}(\S 4.1.1 p 131)
\[
\cZ_K  := (D^2,S^1)^K \in ({\rm Spaces \; with} \; \T^m - {\rm action}) 
\]
\cite{2}(\S 4.1.1 p 131), where now $(X_i,A_i) = (D^2,S^1)$ is a pair with an 
action of the circle.

{\bf Remarks} 

1)  $K \mapsto \cZ_K$ is monoidal, taking joins of simplicial 
complexes to Cartesian products of spaces \cite{2}(\S 4.5.5), and 

2) If $|K| \cong S^{n-1}$ then $\cZ_K$ is a closed topological 
manifold of dimension $n+m$ \cite{2}(\S 4.1.4).

{\bf 1.2.2} Following \cite{2}(\S 3.1.1 p 93), let $\Z[K]$ denote the 
Stanley-Reisner face ring of the simplicial complex $K$: this is the quotient 
of the augmented polynomial algebra $\Z_{1 \leq i \leq m}[v_i] = \Z[v_*]$ 
generated by the vertices of $K$, modulo the ideal generated by relations 
$v_I = \prod_{k \in I} v_k$, where the subset $I \subset K$ is {\bf not} a 
simplex of $K$. Evidently $\Z[K]$ is a $\Z[v_*]$-algebra. 

There has been a great deal of interesting recent work relating combinatorics
(formulated for example in terms of face rings of simplicial complexes) and
the equivariant topology of their associated moment-angle complexes. Some
of this work will be summarized below; \cite{2} contains further references to 
the original literature. Thus \cite{2}(Prop 4.5.5 p 147):
\[
\Tor^*_{\Z[v_*]}(\Z[K],\Z) \cong H^*(\cZ_K;\Z) \;, \; H^*_{\T^m}(\cZ_K;\Z) 
\cong \Z[K]
\]
(as graded rings). Moreover, if $P$ is an associated  omnioriented simple 
polytope as in \S 1.1.1, and we let $T_\Lambda$ denote the kernel of the 
surjection
\[
\Lambda \otimes \R/\Z : \T^m \to \T^n \;,
\]
then $T^n \cong\T^m/\T_\Lambda$ acts on $\cZ_K/T_\Lambda := M$, which carries 
the structure of a complex-oriented quasi-toric manifold \cite{2}(\S 7.3.8 - 
16 p 250); in particular, it is a $\sT$-space. Moreover
\[
H^*_{T^n}(M;\Z) \cong \Z[K]
\]
while $H^*(M;\Z)$ is isomorphic to the quotient of $\Z[K]$ by the ideal 
generated by the image of $\Lambda$ \cite{2}(Th 5.3.1 p 186). These 
isomorphisms have natural interpretations in terms of Eilenberg-Moore spectral 
sequences.

{\bf Remarks} 

1) The composition $T_\Lambda \subset \T^m = \prod_{1 \leq i \leq m}
\T_i \to T_i$ defines a representation $\C(i)$ of $T_\Lambda$ on the $i$th 
circle in the product $\T^m$, and thus a complex line bundle 
\[
L_i := \cZ_K \times_{T_\Lambda} \C(i) \to M
\]
whose Chern classes can naturally be identified with the images of the classes 
$v_i$ in $H^*_{\T^m}(\cZ_P;\Z)$. The images $[v_i]$ of these classes in 
$H^*(M;\Z)$ are those cited above in Guillemin's formula for the K\"ahler 
class of $M$ \cite{2}(\S 7.3.15 p 249).

2) When $K$ is a Gorenstein complex \cite{2}(Th 3.4.2 p 112), $M$ is a 
manifold with fundamental class
\[
\xymatrix{
{H^{2n}_{T_\Lambda}(\cZ_K;\Z) \cong H^{2n}(M;\Z)} \ar[dr] \\
{} & \Z \\
{H^{2n}_{\T^m}(\cZ_K;\Z) \cong \Z[K]} \ar[uu] \ar[ur] \;.}
\]
{\bf 1.2.3} \cite{2}(Th 7.3.15 p 249) We have a $T^n$-equivariant splitting 
\[
T_M \oplus \C^{m-n} \cong \bigoplus_{1 \leq i \leq m} L_i 
\]
of the tangent bundle of $M$ \cite{7},\cite{2}(Th 7.3.15 p 249). This suggests
the interest of the class of complex-oriented manifolds, perhaps 
odd-dimensional, possessing stable splittings of their tangent or normal
bundles as sums of complex line bundles. \bigskip

{\bf 1.3 Toric dynamical systems in theoretical chemistry}

Work on the formal theory of chemical reaction networks dates back at least
to the 1970s \cite{28} but work of Gatermann \cite{21,22,25} and others
around 2000 brought new ideas from toric geometry into the subject. Following 
\cite{15},\cite{47}(\S 2) closely, we summarize fundamental constructions: 

{\bf 1.3.1} A chemical reaction network $\A$ (for alembic) is defined by a 
triple  
\[
\A = \{G, \bk \in \M^n_n(\R), \bY \in \M^s_n(\Z) \} 
\]
where

*) $G$ is a directed graph with a set $V$ of $n$ ordered vertices, called 
reactions or nodes, with edge set $E \subset \{(k,l) \in V \times V \;|\; 
k \neq l \}$ (\ie without closed loops),

*) $\bk = [\kappa^k_l]$ is a matrix of reaction rates satisfying 
$\kappa^k_l > 0$ if $k \neq l$, extended by defining 
\[
\kappa^l_l := - \sum_{k \neq l} \kappa^k_l
\]
(making each row $l$ sum to zero, so $\kappa$ becomes the negative of the 
Laplacian matrix for the weighted directed graph $G$), and 

*) $\bY = [Y^i_k] \in \M^s_n(\Z_{n \geq 0})$ is a non-negative integer matrix, 
with $s$ the number of chemical `species' involved in the network. 

A {\it toric dynamical system} is a chemical reaction network which admits a 
strictly positive constant solution $\bc = (c_1,\dots,c_m)$ of the algebraic 
system $\bk \cdot \bc^\bY = 0$: this is a special kind of steady state of the 
system 
\[
\dot{\bc} = \bY \cdot \bk \cdot \bc^\bY \;,
\]
\ie  
\[
\frac{dc_i(t)}{dt} = \sum_{1 \leq k,l \leq n} Y^i_k \kappa^k_l \cdot 
\prod_{1 \leq j \leq s} c_j^{Y_j^l} = 0 
\]
of ordinary differential equations, a {\it complex balancing} steady state in
the sense of \cite{28}.

Following \cite{48}, a toric ideal in a polynomial ring is a prime ideal
generated by differences of monomials. Theorem 7 of \cite{15} defines a 
variety $V(\A) := V(M_{G,Y})$ associated to a toric ideal in $\Q[\kappa^i_j]$ 
generated by the reaction rates of $\A$, and shows that points $\kappa \in 
V_{>0}(\A)$ characterize toric dynamical systems, whose complex balancing 
states are the positive points of a variety $V(T_{G,Y})$ defined by a toric 
ideal in $\Q[c,\kappa]$. An outline of this construction follows:

{\bf 1.3.2}  Suppose at first that $G$ is strongly connected (\ie there 
is a directed path between any pair of distinct vertices); the matrix 
$\bk$ then has rank $n-1$. Let $(-1)^{n-1} K_i$ denote the determinant of 
a minor of $\bk$ defined by deleting the $i$th row and any column - by 
Stanley's matrix-tree theorem \cite{47}(\S 2.2.1), it doesn't matter 
which - and it follows \cite{15}(Lemma 5) that $\Q[K_1,\dots,K_n] 
\subset \Q[\kappa_*]$ is a polynomial subalgebra. When $G = \coprod_{1 \leq q 
\leq r} G(r)$ is the disjoint union of strongly connected subgraphs, with 
$G(q)$ having $n(q)$ ordered nodes (so $n = \sum n(q)$), the existence of
the polynomial algebra $\Q[K_*]$ generalizes. 

The {\it Cayley matrix} \cite{15}(p 7),\cite{46}
\[
\Cay_G(Y) := \left[\begin{array}{c}
                         [Y(*)] \\ {\bf 1}(*)
                     \end{array}\right] \in M^n_{r+s}(\Z_{\geq 0})
\]
is defined by the row of submatrices $[Y(q)^k_l] \in M^{n(q)}_s(\Z_{\geq 0})$
of $\bY$ (indexed by nodes of $G(q)$; but note the reversal of rows and
columns) across the top, filled in below by appropriate rows of zeroes 
and ones. The ideal $M_{G,Y}$ of $\Q[K_*]$ generated by the elements
\[
\prod K_i^{u_{0,i}} - \prod K_i^{u_{1,i}} \;,
\]
indexed by vectors $\bu_0,\bu_1 \in \N^n$ such that $\Cay_G(Y) \cdot 
(\bu_0 - \bu_1) = 0$, then \cite{15}(Th 9 p 8),\cite{47}(\S 2.2.1 p 16)
defines the toric variety $V(\A)$. The moment map identifies the interior 
of the convex hull of the columns of $\Cay_G(Y)$ with $V_{>0}(\A) \subset 
\R^n$; its codimensions $\delta$ is the {\it deficiency} of the reaction
network. Thus when $\delta = 0$, the moment polytope is an $n$-simplex. 

Interesting applications to Hopf bifurcation theory are considered in 
\cite{22}.

{\bf 1.3.3} The real positive points on such a complex toric variety can be
regarded as asymptotic limits of states of {\it dissipative} biochemical 
systems. As suggested in \cite{17}, the short term mechanics of living 
organisms are essentially metabolic, but in the long run they are low-energy 
information management processes : the egg is the leading term in the 
asymptotic expansion of the chicken. The complex points of these varieties, 
and their topology, will thus not be assigned any immediate physical or 
biological interpretation. \bigskip

{\bf 1.4 The classical characteristic numbers} 

{\bf 1.4.1} The Thom spectrum 
\[
M\U : S^2M\U(k) \to M\U(k+1)
\]
(for the stable unitary group $\U = \cup_{k \geq 1} \U(k)$) is defined by the
Thom spaces $M\U(k)$ of the canonical $\C^k$-bundles $\xi_k$ over the 
classifying spaces $B\U(k)$. The classifying map
\[
\nu : M \to B\U(N-n),\; N \gg 0
\]
of the stable normal bundle of a complex-oriented $2n$-dimensional 
manifold $M \subset \C^N$ defines a class
\[
[\C^N_+ \to M^\nu \to M\U(N-n)] \in \pi_{2n}M\U := M\U_{2n}
\]
in the complex bordism ring $M\U_*$ (with $MU_{\rm odd} = 0$). The Thom 
isomorphism
\[
H_*(M\U;\Z) \to H_*(B\U;\Z)
\]
(of comodules) identifies Hurewicz's ring monomorphism
\[
\pi_*M\U \to \Hom^{-*}(H^*(M\U),H^*(S^0)) \cong H_*(B\U)
\]
with the characteristic number homomorphism which sends $M$ to the evaluation
\[
c^I(\nu)[M] := (\prod_{1 \leq k \leq r} c_k^{i_k}(\nu))[M]
\]
of a degree $2n$ word\begin{footnote}{where $I = 1^{i_1}2^{i_2} \dots r^{i_r},
\; |I| =  \sum k i_k = n$ is a partition of $n$ with $r$ parts}\end{footnote} 
in the Chern classes of the stable normal bundle of $M$ on the fundamental 
class $[M] \in H_{2n}(M;\Z)$. This map becomes an isomorphism after tensoring 
with $\Q$: rationally, complex cobordism theory and the theory of Chern 
numbers coincide.

The group completion 
\[
\coprod_{k \geq 0} B\U(k) \to \Omega B(\coprod_{k \geq 0} B\U(k)) \simeq \Z 
\times B\U \supset 0 \times B\U
\]
defines an $H$-space structure on the classifying space for stable complex
vector bundles. Following Borel and Hirzebruch, we identify its cohomology
\[
H^*(B\U;\Z) \cong \Z_{k \geq 1}[e_k] := \sS^* \subset \Z_{k \geq 1}[x_k] \cong 
\lim_{r \to \infty} H^*(B(T^r);\Z) 
\]
$(|x_k| = 2, \; |e_k| = 2k)$, with the graded polynomial algebra of symmetric
functions \cite{35}(I \S 2.7) generated by elementary symmetric functions, 
identifying $e_k$ with the Chern class $c_k$. In terms of generating functions,
\[
E(t) := \prod_{i \geq 1}(1 + x_i t) = \sum_{j \geq 0} e_jt^j \;,
\]
and the {\it complete} symmetric functions $h_j,\; H(t) = \sum_{j \geq 0} 
h_jt^j$ are defined by $H(t) = E(-t)^{-1}$; thus, for example, $e_0 = h_0 
= 1, \; e_1 = h_1$. The generating function 
\[
tH'(t)H(t)^{-1} = \sum_{k \geq 1} p_k t^k 
\]
defines the power sums $p_k = \sum_{i \geq 1} x_i^k$. 

The ring of symmetric functions has a canonical pair of involutions
generating an action of the Klein four-group: the first is defined by 
$x_i \mapsto -x_i$, sending $e_k$ to $(-1)^k e_k$, while the second maps 
$e_k$ to $- h_k$ (so $p_k \mapsto - p_k$). These correspond to the 
operations which send a stable complex vector $V$ respectively to its 
complex conjugate $V^*$ and to its $H$-space inverse $-V$. It follows 
that the $i$th Chern class of the stable normal bundle $\nu$ of a stably 
almost-complex manifold corresponds to $-h_i$ of its stably 
almost-complex tangent bundle $T_M$ \cite{13}.

The algebra of symmetric functions is moreover self-dual with respect
to the Hall inner product \cite{35}(I \S 4.5), with an isomorphism
\[
H_*(B\U,\Z) \cong \sS_*
\]
defined by sending a word $h_I$ in the complete symmetric functions
$h_i$ to the word $m_I$ written in terms of the corresponding  monomial
symmetric functions $m_i$; the Schur symmetric functions form an an 
orthonormal basis. Note finally that the polynomial algebra $H^*(B\U(1),
\Z) = H^*(\C P^\infty,\Z)$ is dual to the Hopf algebra $H_*(\C P^\infty,
\Z) \cong \Z_{k \geq 1}[b_{k}]$ of divided powers $b_{(k)} =  b^k/k! 
\in \Q[b]$.

{\bf 1.4.2} The Landweber-Novikov Hopf algebra $S_* = \Z_{i \geq 1}[t_i]$
of formal diffeomorphisms of the line at the origin (with $|t_i| = 2i$ and
$t(T) = \sum_{i\geq 0} t_i T^{i+1}$, and an antipode given by power 
series inversion) is defined by the coproduct
\[
(\Delta_S t)(T) = (t \otimes 1)((1 \otimes t)(T)) \in (S \otimes S)_*[[T]]
\]
($|T| = -2$);  Mi\v{s}\v{c}enko's logarithm
\[
\log_{M\U}(T) \; = \; \sum_{n \geq 1} \frac{\C P_{n-1}}{n} T^n \in
M\U^*_\Q[[T]]
\]
for the formal group law on $M\U^*(\C P^\infty) \cong M\U^*[[c]]$ defines 
a coaction
\[
\psi_{M\U}(\log_{M\U}(T)) = \log_{M\U}(t(T))
\]
of $S_*$ on $M\U_*$, yielding an isomorphism
\[
(M\U_*,M\U_* \otimes S_*) \cong (M\U_*,M\U_*M\U)
\]
of Hopf algebroids \cite{44}. The closely related Faa di Bruno Hopf 
algebra 
\[
S_* \otimes \Q = \Q_{k \geq 1}[t_{(k)}] 
\]
is presented in terms of the formal derivatives $t_{(k)} = k!t_{k-1}$ of 
$t(T)$ at 0; this is a confused issue in the literature. Extending either 
by adjoining an invertible degree zero element $t_0^{\pm 1}$ can be used 
to encode the grading. 

The generators $b_i \in M\U_*(\C P^\infty)$ (Kronecker dual to $c^i \in 
M\U^*[[c]]$) define a formal series
\[
b(T) = \sum_{i \geq 0} b_i T^i = \exp(b \log_{M\U}(T)) \in (M\U^\Q_*
\C P^\infty)[[T]]
\]
of degree zero \cite{38}, with coproduct $\Delta b(T) = b(T) 
\otimes_{M\U} b(T)$, making $M\U_*(\C P^\infty)$ an $M\U_*M\U$-comodule 
coalgebra with
\[
\psi_{M\U}b(T) = b(t(T)) \;.
\]

{\bf 1.4.3} The Hurewicz homomorphism 
\[
\B_* \cong M\U_*B\T \to H_*(B\U;\Z) \otimes H_*(B\T;\Z) \cong \sS_* 
\otimes \Z_{i \geq 1}[b_{(i)}]
\]
maps VL Ginzburg's cobordism ring (of integral symplectic manifolds under
non-degenerate cobordism (as in \S 1.1.2 \cite{38}) injectively, sending 
a $2n$-dimensional manifold $(V,\omega)$ to the system 
\[
I \mapsto (-1)^{\sum i_k}(\prod h_k^{i_k}(T^\C_V) \cdot [\omega]^{n - i})
[V] \in \Z 
\]
of characteristic number homomorphisms indexed by unordered partitions 
$I$ of $i \leq n$. Relaxing the integrality constraint on the form 
$[\omega]$ in effect embeds $\B_*$ in $\B_* \otimes \R$. 

Quillen's conventions \cite{43} associate to a compact complex 
$2n$-dimensional $G$-manifold, a class in $M\U^{-2n}(BG)$; similarly, an 
integrally symplectic $2n$-dimensional $\sT$-manifold defines an element 
of $\B^{-2n}(B\T^n)$, with characteristic numbers in 
\[
H^{-2n}(B\T^n;H_*(M\U \wedge B\T)) \;,
\]
generalizing the characteristic numbers for formal Hamiltonian cobordism 
defined in \cite{1}(H.2 p 304). The group homomorphism
\[
\T^n \times \T^n \to \T^n \to \U(n)
\]
defined by multiplication of abelian groups induces a map from $B\T^n 
\times B\T^n$ to $B\U(n)$, and pulling back (the stable inverse of) the 
universal bundle on $B\U$ defines a morphism 
\[
B\T^n_+ \wedge M\U(1)^{\wedge n} \to M\U(n)
\]
of spectra, adjoint to a morphism 
\[
M\U(1)^{\wedge n} \to \Maps(B\T^n_+,M\U(n)) \;.
\]
Taking limits as $n \to \infty$ defines (see \@.1.2)  a morphism from the 
homology xgroups of $M\xi$ to those of $\Maps(B\T^n_+,M\U)$, relating formal 
Hamiltonian characteristic classes to noncommutative symmetric functions. 
\bigskip

{\bf \S II The Baker-Richter spectrum}

{\bf 2.1 New kinds of characteristic classes and numbers}

{\bf 2.1.1} The $2^{n-1}$ {\it ordered} partitions
\[
\bn := n_1 + \cdots + n_k
\]
of $n = \sum n_i$ into nonempty parts define the quasisymmetric function
\[
\la \bn \ra (x_*) := \sum_{0 < n_1 < \dots < n_k} \prod x_{i_j}^{n_j} \in 
\Z[x_*]
\]
(denoted $[n_1,\dots,n_k]$ in \cite{8}, \cf \cite{27}(\S 4)). The subring 
of $\Z[x_*]$ generated by such sums is, by Baker and Richter or 
Hazewinkel's proof of Ditters' conjecture, an evenly graded 
(`quasi-shuffle') Hopf algebra $\sf{QSymm}^*$ over $\Z$, generated over 
$\Q$ by certain Lyndon partitions \cite{8}(\S 2.2); for brevity it will 
here be denoted $\QS^*$. The dual Hopf algebra of {\bf non}commutative 
symmetric functions
\[
\NS_* = {\sf NSymm}_* = \Z_{i \geq 1}\la Z_i \ra
\]
is free associative, with coproduct $\Delta Z_i = \sum_{i=j+k} Z_j \otimes
Z_k$; the word $Z_\bn := \prod_{1 \leq i \leq k} Z_{n_i} \in \NS_{2n}$ is 
dual to $\la \bn \ra \in \QS^{2n}$. We will identify the generators $Z_i$ 
with Cartier's elements $\Lambda_i$, \cf \cite{12}(\S 4.1F eq 155). 
Abelianization defines dual maps
\[
\NS_* \to \sS_*, \; \; \sS^* \to \QS^*
\]
of Hopf algebras\begin{footnote}{Convention: $- \otimes_\Z \Q$ sends the 
graded module $M_* = M^{-*}$ to $M^\Q_* = M^{-*}_\Q$.}\end{footnote}; the 
second sends $h_k \mapsto \la {\boldsymbol{1^k}} \ra$; see further 
\cite{35}(\S 2).]

IM James' construction provides a stable splitting 
\[
\Omega \Sigma B\U(1) \sim \bigvee_{n \geq 0} B\U(1)_+^{\wedge n}
\]
and thus Hopf algebra isomorphisms 
\[
H_*(\Omega \Sigma \C P^\infty;\Z) \cong \NS_*, \; \; H^*(\Omega \Sigma \C 
P^\infty;\Z) \cong \QS^* \;.
\]
{\bf 2.1.2} Regarding a complex line bundle as a complex vector bundle 
defines a map from $B\U(1)$ to $B\U$; composing with the stable inverse 
map of \S 1.4.1 above defines a morphism
\[
\xymatrix{
\Omega \Sigma B\U(1) \ar[r] & B\U \ar[r]^{V \mapsto -V} & B\U }
\]
of loop-spaces. Pulling the canonical stable bundle over $B\U$ back 
defines the $A_\infty$ spectrum $M\xi$ of \cite{9}. A class in 
$\pi_{2n}M\xi$ can then be interpreted as the cobordism class of a 
complex-oriented $2n$-manifold $M$, together with a preferred isomorphism
\[
T_M \simeq \oplus L_i 
\]
of its stable tangent bundle as a sum of complex line bundles; 
omnioriented manifolds provide examples. Forgetting the splitting defines
an abelianization homomorphism $M\xi_* \to M\U_*$; note that expressing 
this in terms of characteristic numbers will involve the second 
involution mentioned above. 

The Hurewicz (ring) homomorphism
\[
M\xi_* := \pi_*M\xi \to H_*(M\xi;\Z) \cong H_*(\Omega \Sigma \C P^\infty;
\Z) \cong \NS_*
\]
is injective, and becomes an isomorphism after tensoring with $\Q$ 
\cite{9},\cite{27}(\S 2.4); it sends $M$ to the homomorphism defined by
\[
(\la \bn \ra([v_*]))[M] := \sum_{0 < n_1 < \dots < n_k}(\prod c_1
(L_{i_j})^{n_j})[M] : \QS^{2n} \to \Z \;,
\]
taking the quasisymmetric function $\la \bn \ra$ (of the Chern classes 
$[v_i]$) to its evaluation on the fundamental class of $M$. 

Similarly, a symplectic $2n$-manifold $(V,\omega)$ with a compatible 
quasitoric structure (for example, a Hamiltonian $\sT$-manifold as in \S 
1.1.2) has characteristic numbers
\[
(V,\omega) \mapsto [\bi \mapsto (\la \bi(T_V) \ra [\omega]^{n-i})[V] 
\cdot b_{(i)}] \in \Hom(\QS^{2n},\R[b]) \cong \NS_{2n}\otimes \R[b] 
\]
indexed by partitions of $i \leq n$. If $(V,\omega)$ is projective toric,
Guillemin's theorem [\S 1.1.3] expresses this invariant in terms of the 
face\ parameters $\lambda_i$. 

The morphism
\[
M\xi \wedge B\T \to M\U \wedge B\T
\]
underlying this map has interesting arithmetic properties \cite{12}(\S4.4)
over $\R$. 
 
{\bf 2.1.3 Remarks:} 

1) A better understanding of toric varieties as limits of toric
manifolds would be very useful. A toric variety associated to a lattice
polytope in $\pi_1(T) \subset (\Lie T)^*$ is smooth \cite{2}(5.2 p 184),
\cite{16}(\S7.1) if the primitive lattice vectors along the edges meeting
at a given vertex form a basis for the lattice. Alexeev \cite{3}(\S 2) has
constructed a moduli space for suitably polarized smooth toric varieties,
but the polarization data restricts their allowed degenerations, making it
difficult to assess the degree to which smooth toric varieties are 
generic. As noted above, the theory of integral characteristic numbers 
of toric manifolds extends nicely to a theory of rational characteristic 
numbers for more general toric varieties. 

2) By the construction in \S 1.2.2, we can associate to quasisymmetric 
functions, characteristic classes in the face rings of simplicial spheres,
and evaluate them to define characteristic numbers for such spheres, 
which are essentially the same as the characteristic numbers for their 
associated (toric) moment-angle manifolds. Since joins map to products 
under that correspondence, this defines a homomorphism from the 
noncommutative graded ring (with respect to joins) of isomorphism classes 
of simplicial spheres, to $\NS_*$. For example, 
\[
\partial \Delta^n \mapsto \sum_{|\bn| = n} Z_\bn \;.
\]
2) The collection $\A$ defining a chemical reaction network, as in \S 1.3,
has an associated projective toric variety $V(\A)$; when smooth, this is a
Hamiltonian $\sT$-manifold, which defines an element in $M\xi_*B\T$ and 
thus an element of $\B_*$. In general it defines an element of $M\xi_*B\T
\otimes \Q$, with image in $\B_* \otimes \R$. Following \S 1.3.2, the 
remark above identifies the class in $M\xi_{2n}$ of a network with 
deficiency $\delta = 0$; see also $\S 1.1.3i$. \bigskip

{\bf 2.2 Formal diffeomorphisms of the noncommutative line}

{\bf 2.2.1} Although composition of formal power series with noncommuting
coefficients is {\bf not} associative, Baker and Richter \cite{9}(Prop 
3.1) show that, with suitable care, $M\xi^*\C P^\infty$ can be understood 
as having a formal one-dimensional Lie group structure defined by a 
coproduct
\[
\Delta_\xi x_\xi = \sum a_{i,j} x_\xi^{i+1} \otimes x_\xi^{j+1}
\]
where, however, $x_\xi \in M\xi^2\C P^\infty$ is {\bf not} a central 
element; note that $M\xi$ possesses a Thom isomorphism for complex vector 
bundles even though it is not an $M\U$-algebra \cite{8}(\S 7).
Nevertheless, their Prop 2.3 uses the Hurewicz homomorphism to define 
an algebra monomorphism
\[
x_\xi \mapsto \Theta(x_\xi) := \sum_{i \geq 0} Z_i x_H^{i+1} = Z(x_H) : 
M\xi^* \C P^\infty \to \NS_*[[x_H]]
\]
with $x_H$ central. Moreover, the rationalization $\Theta \otimes \Q$ is 
an isomorphism, and we have a commutative diagram
\[
\xymatrix{
M\xi^*_\Q\C P^\infty \ar[d]^{\Theta_\Q}_\cong \ar[r]^-{\Delta_\xi} & 
M\xi^*_\Q(\C P^\infty \times \C P^\infty) \cong M\xi^*_\Q\C P^\infty 
\otimes_{M\xi} M\xi^*_\Q\C P^\infty \ar[d]^{\Theta_\Q \otimes 
\Theta_\Q}_\cong \\
\NS_*^\Q[[x_H]] \ar[r]^-{\Delta_\NS} & \NS_*^\Q[[x_H \otimes 1,1 \otimes 
x_H]]} 
\]
defining a formal group law over $\NS_*^\Q := \NS_* \otimes \Q$ with 
central coordinate $x_H$. Similarly, $M\xi^\Q_*\C P^\infty \cong \NS^\Q_*
[b^H_i]_{i \geq 1}$ is a Hopf algebra with binomial coproduct 
\[
b^H(T) = \sum_{i \geq 0} b^H_i T^i, \; \Delta b^H(T) = b^H(T)\otimes_{\NS}
b^H(T) \;,
\]
where $b^H_i$ is Kronecker dual to $x_H^i$.

{\bf Proposition}{{\it The Hurewicz characteristic number construction
defines a commutative diagram
\[
\xymatrix{
M\xi^*\C P^\infty \ar[d] \ar[r] & H^*(\C P^\infty,H_*M\xi) = \NS_*[[x_H]] 
\ar[d] \\
M\U^*\C P^\infty \ar[r] & H^*(\C P^\infty,H_*M\U) = \sS_*[[c]] }
\]
sending $\Theta(x_\xi)$ to $\exp_{M\U}(c)$}.

{\bf 2.2.3} I am indebted to Michiel Hazewinkel for drawing attention to 
work of Brouder, Frabetti, and Krattenthaler on a (neither commutative 
nor cocommutative) Hopf algebra $\FN := (\NS_*,\Delta_\NF)$ of formal 
diffeomorphisms of the noncommutative line, generalizing $(S_*,\Delta_S)$.
It is defined on generators in terms of a formal residue homomorphism
\[
\res_{T=0}(\prod_{i \in \Z} a_i T^i) := a_{-1} \;,
\]
(related to the Stieltjes transform \cite{29}(\S 1.2) used in free probability 
theory) by
\[
\Delta_\NF Z(T) = \sum_{n \geq 0} \res_{U=0} (Z(U) \otimes Z(T)^n U^{-n-1})
\in (\NS \otimes \NS)_*[[T]]
\]
where $T$ and $U$ are central formal indeterminates of degree -2. 
Noncommutative Lagrangian inversion \cite{10} (Th 2.14) defines a 
(nonsymmetric) antipode $\chi_\NF$, \ie $\Delta_\NF \circ \Delta_\NF
\neq {\it id_\NF}$. It is an exercise, using the formal analog of Cauchy's 
theorem, to see that the BFK Hopf algebra abelianizes to the Landweber-Novikov 
algebra. 

It will be helpful to provide a concordance of their notation. They 
work with a complexification $\cHd := \C_{i \geq 1} \langle a_i \rangle$ 
of $\FN$, with $a_i$ corresponding to our $Z_i$, as well as with a 
cocommutative Hopf algebra $\cHi := \C_{i \geq 1}\langle b_i \rangle$ 
(analogous to $\NS_* \otimes \C$, on generators $b_i$ with binomial coproduct 
$\Delta b(T) = b(T) \otimes_\C b(T)$); for example, Cartier's primitive 
elements 
\[
\Psi_k \in \NS_{2k}, \; \Psi(T) := T Z'(T) \cdot Z(T)^{-1}
\]
define such generators by $b(T) := \exp(\Psi(T))$. The algebra homomorphism 
\[
\psi : \cHi \to \cHi \otimes_\C \cHd
\]
defined by
\[
(\psi b)(T) = \sum_{n \geq 0} \res_{U=0} (b(U) \otimes a(T)^n U^{-n-1})
\]
is, by their lemma 3.3, a morphism of coalgebras, defining (a family of)
Hopf algebroid(s) $(\cHi,\cHi \otimes_\C \cHd)$. The choice above defines 
a Hopf algebroid $(\NS^\Q_*, \NS^\Q_* \otimes \FN)$, conjecturally sending 
$\Psi_{k-1}$ to $k^{-1} \C P_{k-1} \in M\U^\Q_{2(k-1)}$. 

{\bf Conjecture} {\it Setting $b^H(T) = \exp(\beta \Psi(T))$ in \S 2.2.1 
defines an isomorphism 
\[
(\NS^\Q_*[\beta],\NS^\Q_*[\beta] \otimes \FN) \to (M\xi^\Q_*\C P^\infty,
M\xi^\Q_*(\C P^\infty \wedge M\xi))
\]
of Hopf algebroids, which specializes at $\beta=1$ to an isomorphism}
\[
(\NS^\Q_*,\NS^\Q_* \otimes \FN) \to (M\xi^\Q_*,M\xi^\Q_*M\xi)
\;. 
\]

{\bf 2.3.1} The cumulant generating function 
\[
\log \mathcal{E}(\exp(tX)) = \sum_{n \geq 1} \kappa_n \frac{t^n}{n!}
\]
of a random variable $X$ in classical statistical mechanics is the 
logarithm of its partition function, \ie of the trace or expectation of 
Boltzmann's factor $\exp(tX)$; up to a factor of $t$ and a normalization, 
this is the Helmholtz free energy of physical chemistry, which 
Schr\"odinger identified with negative entropy \cite{45}(Ch 6 p 74).
Its Legendre transform \cite{36}(\S 4),\cite{37} is roughly Cram\'er's 
risk function, which measures the density of excursions from the mean 
along the normal bundle transverse to an evolutionary trajectory \cite{46}. 

{\bf 2.3.2}  In the theory of free (noncommutative) probability, a state 
(\ie a unital $\C$-valued linear functional $\phi$ on a not necessarily 
commutative algebra) associates to an element $a$ of that algebra, a moment 
distribution $\{\phi(a^n), \; n \geq 0\}$ (with $\phi(a^0) = 1$). Since 
$\NS_*$ is free associative over $\Z$, we can take $a^n$ to be the image 
under a ring homomorphism of Cartier's $\Sigma_n$ (defined in terms of 
generating functins by the relation $\Sigma(T) \cdot Z(-T) = 1$). 

{\bf Definition} Following Voiculescu \cite{41}(Ch XI), let $A$ be a 
commutative ring, with 
\[
g(z) = z^{-1} + \sum_{n \geq 1} g_n z^{-n-1} \in A[[z^{-1}]] \;;
\]
then $g(z) = \gamma(z^{-1})$ defines $\gamma(z) \in A[[z]]$, with formal
inverse $\gamma^{-1}(z) = z + \cdots$ such that 
\[
g(\frac{1}{\gamma^{-1}(z)}) = z \;,
\]
and we can write
\[
K(z) = \sum_{n \geq 0} k_n z^n = \frac{z}{\gamma^{-1}(z)} = 1 + \dots \;.
\]

The free moment generating function for the distribution defined by $a$ is
then 
\[
g(z) = \sum_{n \geq 0} \phi(a^n) z^{-n-1} = z^{-1} + \dots \;,
\]
and its free cumulant generating function $\kappa(z)$ is defined by
$K(z) = 1 + z\kappa(z)$.

Remarkable work \cite{20},\cite{21}(prop 4.2) of Friedrich and McKay (going 
back to \cite{26},\cite{30}(\S 1.3)) interprets the relation between 
one-dimensional formal group laws and their associated Hirzebruch 
multiplicative sequences - in particular 
\[
H(z) = \frac{z}{\log^{-1}_{M\U}(z)} 
\]
in terms of the antipode of the Landweber-Novikov algebra (via Lagrange
inversion) to construct a state on $\NS_* \otimes \C$ with $\log_{M\U}$ as 
its free cumulant. More recent work \cite{29}(Th 3.1), \cite{42}(\S 5.5.4) 
lifts this to a version 
\[
K(x_H) = \chi_\NF(Z(-x_H)) \in \FN[[x_H]] \;.
\]
of a noncommutative cumulant generating function for the distribution 
defined by $\Sigma_n \mapsto a^n$. xs\bigskip

{\bf \S III Afterword} 

In the 1960s V I Arnol'd interpreted Euler's equations for fluid
mechanics on a Riemannian manifold as the geodesic flow on the 
infinite-dimensional Lie group $\SD(M)$ of its volume-preserving 
diffeomorphisms \cite{5,6,18,31}, with associated Hamiltonian systems 
defined by coadjoint orbits as in \S $1.1.3ii$ above: the presence of 
viscosity makes such finite-dimensional analogies reasonable.

On the other hand  D'Arcy Thompson's `simple and most beautiful [inkdrop] 
experiment' (On falling drops, {\it On Growth and Form} p 395-6) 
suggests that the work \cite{4,31,33,40} of Kolmogorov, Arnol'd, Moser and 
others on the stability of quasi-periodic solutions of certain Hamiltonian 
systems (originating in celestial mechanics, but generalized (as above, to 
hydrodynamics, or to chemical reaction networks as in \S 1.3.3)) might have 
applications to the analysis of self-replicating systems in biology.\bigskip

\bibliographystyle{amsplain}

\end{document}